\documentclass[a4paper,10pt]{article}
\usepackage{thesis}
\usepackage{natbib}
\usepackage{graphicx}
\usepackage{transparent}
\usepackage[font={singlespacing,small,it}]{caption}
\begin{document}
\font\myfont=cmr12 at 15pt
\title{{\myfont Travelling wave solutions of the perturbed mKdV equation that represent traffic congestion}}
\date{}
\author{Laura Hattam \footnote{University of Reading}}
\maketitle
\abstract{A well-known optimal velocity (OV) model describes vehicle motion along a single lane road, which reduces to a perturbed modified Korteweg-de Vries (mKdV) equation within the unstable regime. Steady travelling wave solutions to this equation are then derived with a multi-scale perturbation technique. The first order solution in the hierarchy is written in terms of slow and fast variables. At the following order, a system of differential equations are highlighted that govern the slowly evolving properties of the leading solution. Then, it is shown that the critical points of this system signify travelling waves without slow variation. As a result, a family of steady waves with constant amplitude and period are identified. When periodic boundary conditions are satisfied, these solutions' parameters are associated with the driver's sensitivity, $\hat{a}$, which appears in the OV model. For some given $\hat{a}$, solutions of both an upward and downward form exist, with the downward type corresponding to traffic congestion. Numerical simulations are used to validate the asymptotic analysis and also to examine the long-time behaviour of our solutions.}
\section{Introduction}
To minimise congestion it is necessary to understand traffic behaviour, which has led to many traffic related studies with varied perspectives. \citet{nag02} has given an overview of the different methods that analyse vehicle motion. One option is the application of a macroscopic model to characterise vehicle headway. The following optimal velocity (OV) model is an example of this approach,
\beq
\frac{d^2 \Delta x_j}{d t^2}=\hat{a}\left(V(\Delta x_{j+1}(t))-V(\Delta x_{j}(t))-\frac{d \Delta x_j}{d t}\right),\label{traf_hw}
\eeq
where $x_j(t)$ is the position of car $j$ at time $t$, $\Delta x_j=x_{j+1}-x_j$ is car $j$'s headway, $V$ is the car's optimal velocity, $j=0,1,2,\ldots,N$ for $N$ cars on the road and $\hat{a}$ is the driver's sensitivity. This equation was derived by \citet{ban95} to describe vehicle behaviour on a single lane road. As well, they proposed the optimal velocity function,
\beq
V(\Delta x_j(t))=\frac{v_{max}}{2}\left(\tanh(\Delta x_j-h_c)+\tanh(h_c)\right),\label{veldef}\eeq
where $h_c$ is the safety distance and $v_{max}$ is the maximal velocity.

Applying linear stability analysis to (\ref{traf_hw}), a neutral stability line with a critical point is obtained. This line signifies the boundary between two stability regions referred to as metastable and unstable. See \citet{ge05} for further detail.

\citet{mur99} reduced (\ref{traf_hw}) to a perturbed Korteweg-de Vries (KdV) equation within the metastable zone using nonlinear theory. This was the KdV equation with higher order correction terms. They numerically identified traffic solitons propagating over open boundaries, which eventually dissolved. This behaviour is expected within this stability regime since all solutions should tend to the uniform headway. \citet{hat16} studied this problem with periodic boundaries, where cnoidal waves were shown to exist that represented traffic congestion. These solutions were derived with multi-scale analysis and then validated with numerical simulations. Again, these density waves disappeared after some time.

Whereas, solutions corresponding to the unstable region were identified by \citet{kom95}. Beginning with (\ref{traf_hw}) close to the critical point along the neutral stability line, they derived a perturbed modified KdV (mKdV) equation. The leading order solution to this equation was written in terms of Jacobi elliptic functions that were dependent on the elliptic modulation term $m\in[0,1]$. When $m=1$, this solution became the kink soliton, which exhibits the start/stop motion representative of a traffic jam.

\citet{kom95} applied perturbation analysis to seek steady travelling wave solutions of the mKdV traffic model. They established that this solution type only existed when the wave modulus $m$ remained constant and consequently, this wave amplitude and period were fixed. A condition for $m$ as some constant was found in terms of integral constraints that then determined the relationship between $m$ and the wave speed. They referred to the travelling wave solutions with $m=1$ as deformed kink solitons. Otherwise, for constant $m\ne 1$, they were called deformed periodic solitons.

Here, a multi-scale perturbation technique is applied to the perturbed mKdV equation to also identify steady travelling wave solutions. This approach is an adaptation of the method outlined by \citet{hat15} for the steady forced KdV-Burgers equation. Solutions of a similar form to the deformed periodic solitons found by \citet{kom95} are highlighted, which satisfy periodic boundaries. \citet{kom95} proposed that this solution type was always unstable and only deformed kink solitons were observed numerically. The stability of our periodic waves is investigated here.

Such studies as \citet{zhu08} and \citet{zhe12} have numerically examined OV traffic models within the unstable zone, where periodic boundary conditions were imposed. The long-time behaviour was analysed, which revealed solutions that were indicative of mKdV dynamics as kink-like waves appeared. Moreover, \citet{li15} performed numerical simulations over large time intervals of an OV model that described a two-lane system with periodic boundaries. As well, this model was transformed into a perturbed mKdV equation near to the critical point. The numerical results corresponding to this region uncovered steady periodic travelling wave solutions with constant amplitude, mean height and period. Hence, these numerical findings suggest stable periodic solutions to the OV traffic system do propagate within this unstable regime.

The focus of this paper is the derivation of steady travelling wave solutions to (\ref{traf_hw}) and then the analysis of their long-time dynamics. In Section \ref{tfm}, (\ref{traf_hw}) is reduced to a perturbed mKdV equation and then steady travelling wave solutions are determined using a multi-scale perturbation method in Section \ref{pa}. The leading order solution is obtained in terms of Jacobi elliptic functions that depend upon slow and fast variables. At the next order, a dynamical system governing the slow variation of the leading order solution is identified. Then, in Section \ref{fp}, the fixed points of this system are shown to represent a family of steady travelling waves that do not slowly vary. This set of solutions have fixed amplitude, mean height and period. Also, the relations between the solution parameters, the wave speed and the driver's sensitivity are established due to implementing periodic boundary conditions. Lastly, in Section \ref{sps}, the highlighted periodic asymptotic solutions are compared with numerical results.

\section{Traffic Flow Model\label{tfm}}
We now outline how the OV model (\ref{traf_hw}) was transformed by \citet{ge05} into a perturbed mKdV equation. This then becomes the steady perturbed Gardner equation when we seek steady travelling wave solutions.

Firstly, \citet{ge05} deduced for the system (\ref{traf_hw}), the linear stability criteria
\beq
\hat{a}\ge\hat{a}_s={2V^{'}(h)},\label{linstab}\eeq
where $h$ is the uniform headway. When this is met, the steady state of $\Delta x_j(t)=h$ is stable. The curve defined by $\hat{a}=\hat{a}_s$ is the `neutral stability line', which indicates the onset of instability. This curve's critical point occurs at $h=h_c$ and $\hat{a}=\hat{a}_c=2V^{'}(h_c)$. The region neighbouring this point is where the perturbed mKdV equation applies.

Next, \citet{ge05} used the change of variables
\beq
\bar{x}=\epsilon(j+V^{'}(h_c)t),\quad \bar{t}=\epsilon^3 g_1 t,\quad \epsilon^2=(\hat{a}_c/\hat{a})-1,\quad 0<\epsilon\ll 1,\label{jttoxtb}\eeq
and let
\beq
\Delta x_j(t)=h_c+\epsilon\sqrt{\frac{g_1}{g_2}}R,\label{Rdef}\eeq
where they introduced
\beq
g_1=V^{'}(h_c)/6,\quad g_2=-V^{'''}(h_c)/6,\quad g_3=V^{'}(h_c)/2,\quad g_4=V^{'}(h_c)/8,\quad g_5=V^{'''}(h_c)/12.\eeq

As a result, (\ref{traf_hw}) reduced to
\beq
R_{\bar{t}}-R_{\bar{x}\bar{x}\bar{x}}+3R^2 R_{\bar{x}}+\epsilon\left(\frac{g_3}{\sqrt{g_1}}R_{\bar{x}\bar{x}}+\frac{g_4}{\sqrt{g_1}}R_{\bar{x}\bar{x}\bar{x}\bar{x}}+\sqrt{g_1}\frac{g_5}{g_2}\partial_{\bar{x}\bar{x}}(R^3)\right)=0,\label{Re}
\eeq
where $O(\epsilon^2)$ terms are ignored, which is a perturbed mKdV equation. Since $\epsilon>0$ is small and therefore $\hat{a}<\hat{a}_c$, the stability criteria (\ref{linstab}) does not hold and the solutions to (\ref{Re}) are unstable. Thus, this is the unstable regime.

Here, steady travelling wave solutions of (\ref{Re}) are sought. To identify this solution type, we set
\bdis
\tilde{x}=\bar{x}-\omega \bar{t},\edis
where $\omega$ is the constant wave speed. As well, if $\sqrt{\omega}u=R+\sqrt{\frac{\omega}{3}}$, (\ref{Re}) becomes
\beq
\lambda u_{\tilde{x}\tilde{x}\tilde{x}}-\gamma u^2u_{\tilde{x}}+\nu u u_{\tilde{x}}+\epsilon G(u,\tilde{x})=0,\label{pmkdv}
\eeq
where
\beq
\lambda=1,\quad \gamma=3\omega,\quad \nu=2\sqrt{3}\omega,\label{eqnpar}\eeq
and
\beq
G(u,\tilde{x})=-\left(\frac{g_3}{\sqrt{g_1}}u_{\tilde{x}\tilde{x}}+\frac{g_4}{\sqrt{g_1}}u_{\tilde{x}\tilde{x}\tilde{x}\tilde{x}}+\omega\sqrt{g_1}\frac{g_5}{g_2}\partial_{\tilde{x}\tilde{x}}\left(\left(u-\frac{1}{\sqrt{3}}\right)^3\right)\right).
\label{Gdef}\eeq
The system (\ref{pmkdv}) is the steady perturbed Gardner equation. The parameters $\lambda,\gamma,\nu$ have been introduced so that initially the analysis is generalised.
\section{Perturbation Analysis\label{pa}}
The modulation theory detailed by \citet{hat15} for the steady perturbed KdV equation is now applied to (\ref{pmkdv}). However, since an additional cubic nonlinear term must be considered, the modulation theory for the Gardner equation is also used, which was outlined by \citet{kam12}. As a result of this perturbation analysis, a leading order solution to (\ref{pmkdv}) is highlighted that varies with slow and fast variables. Then at the next order, equations are found that describe the slow evolution of this solution.

To begin, let
\beq
X=\epsilon \tilde{x},\quad \theta_X=\frac{1}{\epsilon}k(X),\eeq
and
\beq
u(\tilde{x})=u_0(\theta,X)+\epsilon u_1(\theta,X) +\epsilon^2 u_2(\theta,X)+\ldots,\eeq
where $X$ and $\theta$ are `slow' and `fast' variables respectively.

Consequently, (\ref{pmkdv}) takes the form at first and second order
\bsub
\label{12orde}
\baln
&O(1): \lambda k^2 u_{0,\theta\theta\theta} -\gamma u_0^2 u_{0,\theta}+\nu u_0u_{0,\theta}=0,\label{12ordea}\\
&O(\epsilon): \lambda k^3 u_{1,\theta\theta\theta}-\gamma k\left(u_0^2u_1\right)_{\theta}+\nu k(u_0u_1)_{\theta}+g=0,\label{12ordeb}
\end{align}
\esub
where
\beq
g(\theta,X)=G(\theta,X)+3\lambda k^2 u_{0,\theta\theta X}-\gamma u_0^2u_{0,X}+\nu u_0u_{0,X}+3\lambda k k_X u_{0,\theta\theta}.\label{gdef}\eeq

The integration of (\ref{12ordea}) twice gives
\beq
\frac{6\lambda k^2}{\gamma}u_{0,\theta}^2=u_0^4-\frac{2\nu}{\gamma}u_0^3+\frac{12}{\gamma}\hat{C}u_0+\frac{12}{\gamma}\hat{D},\label{CDeqn}\eeq
where $\hat{C}$ and $\hat{D}$ are integration constants. Next, let
\beq
Q(u_0)=u_0^4-\frac{2\nu}{\gamma}u_0^3+\frac{12}{\gamma}\hat{C}u_0+\frac{12}{\gamma}\hat{D},\label{qdef}\eeq
where $Q$ is a polynomial of order $4$. Suppose that $a\le b\le c\le d$ are the roots of this polynomial, then
\beq
u_0^4-\frac{2\nu}{\gamma}u_0^3+\frac{12}{\gamma}\hat{C}u_0+\frac{12}{\gamma}\hat{D}=Q(u_0)=(u_0-a)(u_0-b)(u_0-c)(u_0-d).\label{Qeqn}\eeq
By expanding the righthand side of (\ref{Qeqn}) and equating like terms, we find
\bsub
\label{pardef}
\baln
&\frac{2\nu}{\gamma}=a+b+c+d,\\
&\frac{12 \hat{C}}{\gamma}=-(acd+abd+abc+bcd),\\
&\frac{12\hat{D}}{\gamma}=abcd,\\
&0=bc+ac+cd+ab+bd+ad.
\end{align}\esub
So the equation parameters are dependent upon the roots $a,b,c,d$. As well, since $Q(u_0=a,b,c,d)=0$, then
\bsub
\label{abcdeqn}
\baln
&a\hat{C}+\hat{D}=\frac{\nu a^3}{6}-\frac{\gamma a^4}{12},\label{abcdeqna}\\
&b\hat{C}+\hat{D}=\frac{\nu b^3}{6}-\frac{\gamma b^4}{12},\label{abcdeqnb}\\
&c\hat{C}+\hat{D}=\frac{\nu c^3}{6}-\frac{\gamma c^4}{12},\label{abcdeqnc}\\
&d\hat{C}+\hat{D}=\frac{\nu d^3}{6}-\frac{\gamma d^4}{12}.\label{abcdeqnd}
\end{align}\esub

The solution to (\ref{12ordea}) can also be written in terms of these roots, such that
\beq
u_0(\theta,X)=\frac{c e +d\operatorname{sn}^2\left(\beta(\theta-\theta_0);m\right)}{e+\operatorname{sn}^2\left(\beta(\theta-\theta_0);m\right)},\label{u0soln}
\eeq
where
\bsub
\label{solpar}
\baln
&\frac{24\lambda \beta^2k^2}{\gamma}=(a-c)(b-d),\label{solpara}\\
&e=-\left(\frac{b-d}{b-c}\right),\label{solparb}\\
&m^2=\frac{(a-d)(b-c)}{(a-c)(b-d)},\label{solparc}\\
&\beta=K(m)/P.\label{solpard}
\end{align}\esub
The period of (\ref{u0soln}) is $2P$, where $P$ is some fixed constant. The function $\operatorname{sn}$ is the Jacobi elliptic function, $m\in[0,1]$ is its elliptic modulus and the function $K(m)$ is the complete elliptic integral of the first kind. The parameters $a,b,c,d,m,\theta_0$ are all dependent on the slow variable $X$.

To ensure the next order solution, $u_1$, has the period $2P$, the periodicity conditions
\bsub
\label{intcond}
\baln
&\frac{1}{2P}\int_{\theta_1}^{\theta_2}g(\theta,X)d\theta=0,\label{intconda}\\
&\frac{1}{2P}\int_{\theta_1}^{\theta_2}g(\theta,X)u_0 d\theta=0,\label{intcondb}
\end{align}
\esub
are imposed, where $\theta_1=-P+\theta_0$ and $\theta_2=P+\theta_0$. Note that $u_0$ and its derivatives with respect to $\theta$ are assumed to be periodic over the domain $\theta\in[\theta_1,\theta_2]$. The integral (\ref{intcondb}) written in full is
\beq
\frac{1}{2P}\int_{\theta_1}^{\theta_2}G u_0d\theta+\frac{1}{2P}I_b=0,\label{intb}\eeq
where
\beq
\bspl
I_b&=\int_{\theta_1}^{\theta_2}\left(3\lambda( k^2 u_{0,\theta\theta X}+k k_X u_{0,\theta\theta})-\gamma u_0^2u_{0,X}+\nu u_0u_{0,X}\right) u_0 d\theta\\
&=\partial_X\int_{\theta_1}^{\theta_2}\left(-\frac{3\lambda k^2}{2}u_{0,\theta}^2-\frac{\gamma}{4}u_0^4+\frac{\nu}{3} u_0^3\right) d\theta.\\
\end{split}\label{Ib}
\eeq

Next, by manipulating (\ref{CDeqn}), it can be shown
\bsub
\label{specints}
\baln
&\frac{3\gamma}{2P}\int_{\theta_1}^{\theta_2}u_0^4d\theta=\frac{5\nu}{2P}\int_{\theta_1}^{\theta_2}u_0^3d\theta-\frac{18\hat{C}}{2P}\int_{\theta_1}^{\theta_2}u_0d\theta-12\hat{D},\label{specintsa}\\
&\frac{18\lambda k^2}{2P}\int_{\theta_1}^{\theta_2}u_{0,\theta}^2d\theta=-\frac{\nu}{2P}\int_{\theta_1}^{\theta_2}u_0^3 d\theta+\frac{18}{2P}\hat{C}\int_{\theta_1}^{\theta_2}u_0 d\theta+24\hat{D}.\label{specintsb}
\end{align}
\esub
Hence,
\bdis
-\frac{18\lambda k^2}{2P}\int_{\theta_1}^{\theta_2}u_{0,\theta}^2d\theta+\frac{1}{2P}\int_{\theta_1}^{\theta_2}(-3\gamma u_0^4+4\nu u_0^3)d\theta=-12\hat{D},\edis
and therefore, (\ref{Ib}) takes the form
\bdis\bspl
\frac{1}{2P}I_b
&=-\hat{D}_X.
\end{split}\edis

As well, the integral (\ref{intconda}) can be written
\beq
\frac{1}{2P}\int_{\theta_1}^{\theta_2}Gd\theta+\frac{1}{2P}I_a=0,\eeq
where
\bdis\bspl
I_a&=\int_{\theta_1}^{\theta_2}\left(3\lambda (k^2 u_{0,\theta\theta X}+kk_Xu_{0,\theta\theta})-\frac{\gamma}{3}(u_0^3)_X+\frac{\nu}{2}(u_0^2)_X\right)d\theta\\
&=\partial_X\int_{\theta_1}^{\theta_2}\left(-\frac{\gamma}{3}u_0^3+\frac{\nu}{2}u_0^2\right)d\theta,\\
&=2\hat{C}_X,
\end{split}
\edis
using (\ref{CDeqn}).

Thus, the integral conditions (\ref{intcond}) reduce to
\bsub
\label{ChDhe}
\baln
\hat{C}_X=-\frac{1}{2P}\int_{\theta_1}^{\theta_2}G(\theta,X) d\theta,\label{ChDhea}\\
\hat{D}_X=\frac{1}{2P}\int_{\theta_1}^{\theta_2}G(\theta,X) u_0 d\theta.\label{ChDheb}
\end{align}
\esub
This system determines the slow variation of the leading order solution (\ref{u0soln}) for some given function $G(\theta,X)$.

\subsection{Application to the traffic flow model\label{trafapp}}
So to locate steady travelling wave solutions of the traffic model (\ref{traf_hw}), $G$ is now defined by (\ref{Gdef}). Firstly,
\bdis\begin{split}
\int_{\theta_1}^{\theta_2} G(\theta,X) d\theta&=-\omega\sqrt{g_1}\frac{g_5}{g_2}k^2\left[3\left(u_0-\frac{1}{\sqrt{3}}\right)^2u_{0,\theta}\right]_{\theta_1}^{\theta_2}=0,
\end{split}\edis
where $O(\epsilon)$ terms are ignored. Therefore, if $\kappa_1$ is some constant, from (\ref{ChDhea}),
\beq
\hat{C}=\kappa_1.
\eeq
Next let us set $\hat{D}_X=\tilde{I},$ so that
\bdis
\tilde{I}=\frac{1}{2P}\int_{\theta_1}^{\theta_2} G(\theta,X) u_0 d\theta.\edis
Given (\ref{Gdef}), this integral becomes
\beq\bspl
\tilde{I}
&=-\frac{1}{2P}\int_{\theta_1}^{\theta_2}u_0\left(\frac{g_3}{\sqrt{g_1}}k^2u_{0,\theta\theta}+\frac{g_4}{\sqrt{g_1}}k^4u_{0,\theta\theta\theta\theta}+\omega\sqrt{g_1}\frac{g_5}{g_2}k^2\partial_{\theta\theta}\left(\left(u_0-\frac{1}{\sqrt{3}}\right)^3\right)\right)d\theta\\
&=-\frac{1}{2P}\int_{\theta_1}^{\theta_2}\left(\frac{g_3}{\sqrt{g_1}}k^2u_0u_{0,\theta\theta}+\frac{g_4}{\sqrt{g_1}}k^4u_0u_{0,\theta\theta\theta\theta}-\omega\sqrt{g_1}\frac{g_5}{g_2}k^2u_{0,\theta}\partial_{\theta}\left(\left(u_0-\frac{1}{\sqrt{3}}\right)^3\right)\right)d\theta\\
&=-\frac{1}{2P}\int_{\theta_1}^{\theta_2}\left(-\frac{g_3}{\sqrt{g_1}}k^2u_{0,\theta}^2+\frac{g_4}{\sqrt{g_1}}k^4u_0u_{0,\theta\theta\theta\theta}-3\omega\sqrt{g_1}\frac{g_5}{g_2}k^2u_{0,\theta}^2\left(u_0^2-\frac{2}{\sqrt{3}}u_0+\frac{1}{3}\right)\right)d\theta.
\end{split}\label{tilI}
\eeq
Hence, to write $\tilde{I}$ in full we must solve the integrals
\bdis
\tilde{I}_1=\frac{1}{2P}\int_{\theta_1}^{\theta_2}u_{0,\theta}^2d\theta,\quad
\tilde{I}_2=\frac{1}{2P}\int_{\theta_1}^{\theta_2}u_0u_{0,\theta}^2d\theta,\quad
\tilde{I}_3=\frac{1}{2P}\int_{\theta_1}^{\theta_2}u_0^2u_{0,\theta}^2d\theta,\quad
\tilde{I}_4=\frac{1}{2P}\int_{\theta_1}^{\theta_2}u_0u_{0,\theta\theta\theta\theta}d\theta.
\edis

As a result of manipulating (\ref{CDeqn}), the integrals $\tilde{I}_1,\;\tilde{I}_2$ and $\tilde{I}_3$ are determined. Omitting the details here, we arrive at
\bsub
\label{tilI123}
\baln
&\tilde{I}_1=\frac{1}{\lambda k^2}\left(-\frac{\nu^2}{12\gamma}\alpha_2+\left(\alpha_1+\frac{\nu}{6\gamma}\right)\hat{C}+\frac{4}{3}\hat{D}\right),\label{tilI123a}\\
&\tilde{I}_2=\frac{\gamma}{6\lambda k^2}\left(\alpha_2\left(\frac{9\hat{C}}{2\gamma}-\frac{5\nu^3}{8\gamma^3}\right)
+\alpha_1\left(\frac{6\hat{D}}{\gamma}+\frac{3\nu\hat{C}}{2\gamma^2}\right)+\frac{5\nu^2\hat{C}}{4\gamma^3}+\frac{\nu\hat{D}}{\gamma^2}\right),\label{tilI123b}\\
&\tilde{I}_3=\frac{\gamma}{6\lambda k^2}\left(\alpha_2\left(\frac{69\hat{C}\nu}{10\gamma^2}-\frac{7\nu^4}{8\gamma^4}+\frac{24\hat{D}}{5\gamma}\right)
+\alpha_1\left(\frac{6\nu\hat{D}}{5\gamma^2}+\frac{21\nu^2\hat{C}}{10\gamma^3}\right)+\frac{35\nu^3\hat{C}}{20\gamma^4}+\frac{7\nu^2\hat{D}}{5\gamma^3}-\frac{54\hat{C}^2}{5\gamma^2}\right),\label{tilI123c}
\end{align}
\esub
where
\beq
\alpha_1(m,c,d,e)=\frac{1}{2P}\int_{\theta_1}^{\theta_2}u_0 d\theta,\quad \alpha_2(m,c,d,e)=\frac{1}{2P}\int_{\theta_1}^{\theta_2}u_0^2 d\theta.
\eeq
Refer to the Appendix for the evaluation of these integrals.

Next, the integral $\tilde{I}_4$ can be written, from (\ref{12ordea}),
\bdis\bspl
\tilde{I}_4=\frac{1}{2P}\int_{\theta_1}^{\theta_2}u_0u_{0,\theta\theta\theta\theta}d\theta&=-\frac{1}{2P}\int_{\theta_1}^{\theta_2}u_{0,\theta}u_{0,\theta\theta\theta}d\theta\\
&=-\frac{1}{2P}\int_{\theta_1}^{\theta_2}u_{0,\theta}\left(\frac{\gamma}{\lambda k^2}u_0^2u_{0,\theta}-\frac{\nu}{\lambda k^2}u_0u_{0,\theta}\right)d\theta\\
&=-\frac{\gamma}{\lambda k^2}\tilde{I}_3+\frac{\nu}{\lambda k^2}\tilde{I}_2.
\end{split}
\edis

Thus, (\ref{tilI}) reduces to
\beq
\tilde{I}=k^2\tilde{I}_1\left(\frac{g_3}{\sqrt{g_1}}+\frac{\sqrt{g_1}g_5 \omega}{g_2}\right)
+k^2\tilde{I}_2\left(-\frac{\nu g_4}{\lambda \sqrt{g_1}}-\frac{2\sqrt{3}\omega\sqrt{g_1}g_5}{g_2}\right)
+k^2\tilde{I}_3\left(\frac{g_4\gamma}{\sqrt{g_1}\lambda}+\frac{3\omega\sqrt{g_1}g_5}{g_2}\right),\label{tilI2nd}
\eeq
where $\tilde{I}_1$, $\tilde{I}_2$ and $\tilde{I}_3$ are defined using (\ref{tilI123}).

\section{Fixed Points\label{fp}}
In Section \ref{pa}, the leading order solution (\ref{u0soln}) to the traffic model (\ref{traf_hw}) was obtained. Moreover, its slow variation was found to be governed by the system
\bsub
\label{sysfp}
\baln
&\hat{C}_X=0,\\
&\hat{D}_X=\tilde{I},
\end{align}
\esub
where $\tilde{I}$ is given by (\ref{tilI2nd}). The critical points of this system occur when $\hat{D}_X=\tilde{I}=0$, so that $\hat{D}=\kappa_2$, where $\kappa_2$ is some constant. We now want to highlight these fixed point solutions and show that they have constant wave amplitude, mean height and period.

Firstly, taking the derivative of (\ref{abcdeqna}) with respect to $X$ gives
\bdis
\hat{D}_X=\frac{1}{6}\left(3a^2\nu-2a^3\gamma-6\kappa_1\right)a_X,\edis
using $\hat{C}=\kappa_1$. Similarly, (\ref{abcdeqnb})-(\ref{abcdeqnd}) take the form
\bdis
\hat{D}_X=\frac{1}{6}\left(3b^2\nu-2b^3\gamma-6\kappa_1\right)b_X=\frac{1}{6}\left(3c^2\nu-2c^3\gamma-6\kappa_1\right)c_X
=\frac{1}{6}\left(3d^2\nu-2d^3\gamma-6\kappa_1\right)d_X.\edis
So, if $\hat{D}_X=0$, then
\beq
a_X=b_X=c_X=d_X=0,\label{confp}\eeq
(assuming $3a^2\nu-2a^3\gamma-6\kappa\neq0$ etc.). As well, from (\ref{solparc}) it can be shown
\bdis\bspl
\frac{2}{m}m_X=&a_X\left(\frac{1}{a-d}-\frac{1}{a-c}\right)+b_X\left(\frac{1}{b-c}-\frac{1}{b-d}\right)\\
&+c_X\left(\frac{1}{a-c}-\frac{1}{b-c}\right)+d_X\left(\frac{1}{b-d}-\frac{1}{a-d}\right).
\end{split}\edis
Therefore, when $\hat{D}=\kappa_2$ ($\tilde{I}=0$), (\ref{confp}) is satisfied, and therefore
\beq
m_X=0.\label{mfixed}\eeq
Hence, the fixed point solutions of the system (\ref{sysfp}) correspond to the leading order solution (\ref{u0soln}) with constant wave modulus $m$ and constant solution parameters $a,b,c,d,\beta$. So then, these particular steady travelling waves do not have slowly varying properties, which means their amplitude, mean height and period remain fixed. Therefore, they are of a similar form to the deformed periodic solitons detailed by \citet{kom95}.

It is necessary to determine when $\tilde{I}=0$ to locate these fixed point solutions. Firstly, if the equation parameters $\lambda$ and $\nu$ are defined using (\ref{eqnpar}) and with $\hat{C}=\kappa_1$, $\hat{D}=\kappa_2$, (\ref{tilI2nd}) is written
\beq\bspl
\tilde{I}
&=k^2\tilde{I}_1\frac{g_3}{\sqrt{g_1}}
+k^2\omega\left(\tilde{I}_1\frac{\sqrt{g_1}g_5 }{g_2}+(3\tilde{I}_3-2\sqrt{3}\tilde{I}_2)\left(\frac{g_4}{\sqrt{g_1}}+\frac{\sqrt{g_1}g_5}{g_2}\right)\right),
\end{split}
\eeq
where
\bsub
\baln
&\tilde{I}_1=\frac{\gamma}{k^2}\left(\alpha_1\frac{\kappa_1}{\gamma}-\frac{\alpha_2}{9}+\frac{\kappa_1}{3\sqrt{3}\gamma}+\frac{4\kappa_2}{3\gamma}\right),\\
&\tilde{I}_2=\frac{\gamma}{6 k^2}\left(\alpha_1\left(\frac{3\kappa_1}{\sqrt{3}\gamma}+\frac{6\kappa_2}{\gamma}\right)
+\alpha_2\left(\frac{9\kappa_1}{2\gamma}-\frac{5}{\sqrt{3}}\right)
+\frac{5\kappa_1}{\sqrt{3}\gamma}+\frac{2\kappa_2}{\sqrt{3}\gamma}\right),\\
&\tilde{I}_3=\frac{\gamma}{6 k^2}\left(\alpha_1\left(\frac{14\kappa_1}{5\gamma}+\frac{12\kappa_2}{5\sqrt{3}\gamma}\right)
+\alpha_2\left(\frac{69\kappa_1}{5\sqrt{3}\gamma}+\frac{24\kappa_2}{5\gamma}-\frac{14}{9}\right)
+\frac{14\kappa_1}{3\sqrt{3}\gamma}+\frac{28\kappa_2}{15\gamma}-\frac{54}{5}\left(\frac{\kappa_1}{\gamma}\right)^2\right).
\end{align}
\esub
Next, by setting $\tilde{I}=0$, a definition for the constant wave speed is obtained as a function of $\kappa_1/\gamma$, $\kappa_2/\gamma$, $\alpha_1(m,c,d,e)$, $\alpha_2(m,c,d,e)$ and $h_c$, such that
\beq
\omega=-\frac{\rho_1}{\rho_2},\label{om_def}
\eeq
where
\bdis
\rho_1=\frac{g_3}{\sqrt{g_1}}\tilde{I}_1,\edis
and
\bdis
\rho_2=\tilde{I}_1\frac{\sqrt{g_1}g_5 }{g_2}+(3\tilde{I}_3-2\sqrt{3}\tilde{I}_2)\left(\frac{g_4}{\sqrt{g_1}}+\frac{\sqrt{g_1}g_5}{g_2}\right).
\edis

As an aside, if the minimum and maximum headway are set to $h_{min}$ and $h_{max}$ respectively, we want solutions with $h_{min}\le h_c\le h_{max}$ ($\Delta x_j=h_c+O(\epsilon)$), which suggests $R_{min}\le 0 \le R_{max}$ (refer to (\ref{Rdef})). Consequently, given $u_0\in[b,c]$, any solution must have $b\le 1/\sqrt{3}$ and $c\ge 1/\sqrt{3}$. However, from (\ref{abcdeqnb}) and (\ref{abcdeqnc}), so that both criterion are satisfied,  it is required
\beq
\frac{\kappa_2}{\gamma}=\frac{1}{36}-\frac{1}{\sqrt{3}}\frac{\kappa_1}{\gamma},\label{Dhatdef}
\eeq
where $\kappa_1/\gamma$ and $\kappa_2/\gamma$ are some constants. A definition for $\kappa_2/\gamma$ as a function of $\kappa_1/\gamma$ has now been found. Next, rewriting (\ref{qdef}) in terms of $\kappa_1/\gamma$ and $\kappa_2/\gamma$ gives
\beq
{Q}(z)=z^4-\frac{4}{\sqrt{3}}z^3+12\frac{\kappa_1}{\gamma}z+12\frac{\kappa_2}{\gamma}.\label{Qdef2}\eeq
Finding the four roots of this polynomial determines the solution parameters $a,b,c,d$, and then $m$ using (\ref{solparc}), for some given $\kappa_1/{\gamma}$ and where $\kappa_2/{\gamma}$ is defined by (\ref{Dhatdef}). So, $a,b,c,d$, as well as the wave modulus $m$ are identified as a function of $\kappa_1/{\gamma}$, which are depicted along the top panel of Figure \ref{Ch_m_plot}. Given this relation between $m,a,b,c,d$ and the parameter $\kappa_1/\gamma$, we can then obtain the wave speed (\ref{om_def}) as a function of $\kappa_1/\gamma$ (with $h_c$ set to $4$ throughout), which is displayed on the bottom left of Figure \ref{Ch_m_plot} (the properties $\alpha_1$ and $\alpha_2$ depend upon $m$, $c$, $d$ and $e=-(b-d)/(b-c)$). It should be noted that the shift $\theta_0$ is arbitrary since it has no influence on the other solution parameters. So we set $\theta_0=0$ for the remainder of the paper.
\subsection{Periodic boundaries}

\begin{figure}
\begin{center}
\includegraphics[width=2.8in]{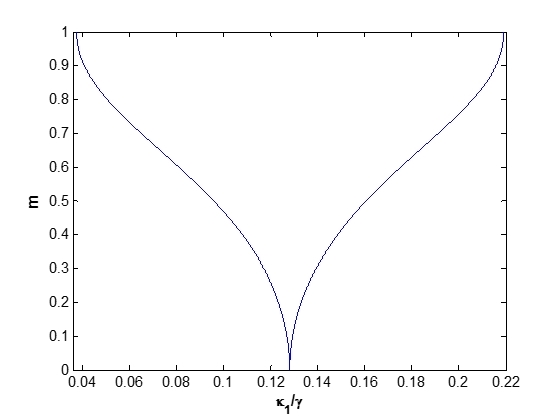}
\includegraphics[width=2.8in]{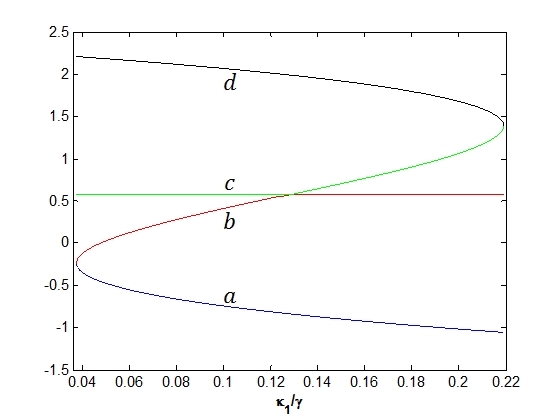}\\
\includegraphics[width=2.8in]{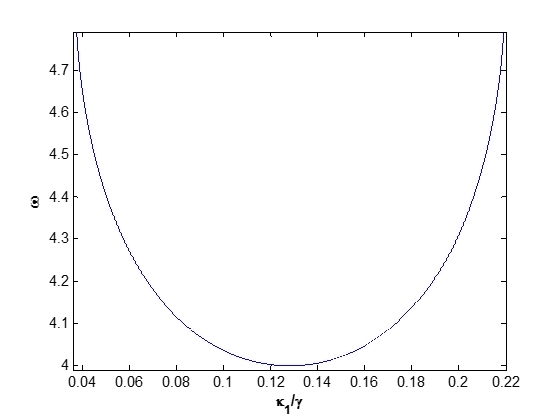}
\includegraphics[width=2.8in]{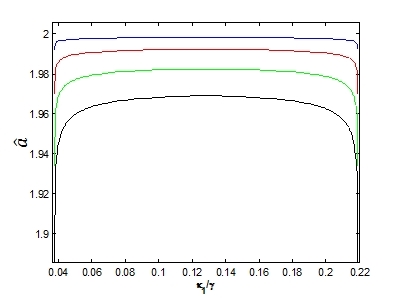}
\caption{The parameter space for the periodic asymptotic solutions when $h_c=4$. Top left: The plot of $m$, (\ref{solparc}), as a function of $\kappa_1/\gamma$. Top right: The plot of $a$ (blue), $b$ (red), $c$ (green), $d$ (black) as a function of $\kappa_1/\gamma$, which are the roots of (\ref{Qdef2}). Bottom left: The plot of $\omega$, (\ref{om_def}), as a function of $\kappa_1/\gamma$. Bottom right: The plot of the driver's sensitivity $\hat{a}$, (\ref{ahdef}), as a function of $\kappa_1/\gamma$, where $N=100$ and blue: $n=1$, red: $n=2$, green: $n=3$, black: $n=4$. Defining $\hat{a}$ with (\ref{ahdef}) ensures periodic boundary conditions are satisfied.\label{Ch_m_plot}}
\end{center}
\end{figure}

The parameters of the leading order solution (\ref{u0soln}) have been defined in terms of the constant $\kappa_1/\gamma$ (see Figure \ref{Ch_m_plot}). By now using the variables $j$ and $t$ that appear in the traffic model (\ref{traf_hw}) and applying periodic boundary conditions, the connection between this constant and the driver's sensitivity, $\hat{a}$, is also established.

To begin, (\ref{u0soln}) is written in terms of $j$ and $t$,
\beq
u_0(j,t)=c+(d-c)\frac{\operatorname{sn}^2\left(\beta\left(k\left(\epsilon(j+V^{'}(h_c)t)-\omega g_1\epsilon^3 t\right)-\theta_0\right);m\right)}
{e+\operatorname{sn}^2\left(\beta\left(k\left(\epsilon(j+V^{'}(h_c)t)-\omega g_1\epsilon^3 t\right)-\theta_0\right);m\right)},\label{u0jt}
\eeq
where (\ref{jttoxtb}) was used. To ensure this solution satisfies periodic boundary conditions, it is necessary for
\beq
\frac{2P}{k\epsilon}n=N,\label{pcon}\eeq
where $n$ is some positive integer representing the number of oscillations over the domain $j\in[0,N]$. However, from (\ref{solpara}), it is known
\beq
\left(\frac{P}{k}\right)^2=\frac{8 K(m)^2}{\omega (a-c)(b-d)}.\label{pcon2}
\eeq
Therefore, combining (\ref{pcon}) and (\ref{pcon2}), as well as using $\epsilon^2=(\hat{a}_c/\hat{a})-1$, we arrive at
\bdis
\frac{\hat{a}_c}{\hat{a}}=\frac{32 K(m)^2n^2}{\omega (a-c)(b-d) N^2}+1.\edis
Rearranging this, a definition for the driver's sensitivity is found,
\beq
\hat{a}=\frac{a_c\omega (a-c)(b-d) N^2}{32 K(m)^2n^2+\omega (a-c)(b-d)N^2}.\label{ahdef}
\eeq
For our fixed point solutions, the relationships between $m,a,b,c,d,\omega$ and $\kappa_1/\gamma$ have been obtained (see Figure \ref{Ch_m_plot}). Using these and (\ref{ahdef}), $\hat{a}$ as a function of $\kappa_1/\gamma$ is determined for some fixed $n/N$. This relation is plotted on the bottom right of Figure \ref{Ch_m_plot} for various $n/N$ values. Note that when (\ref{ahdef}) holds, the solution will satisfy periodic boundary constraints. By relating $\hat{a}$ to the asymptotic analysis then enables us to compare numerical solutions of the OV system (\ref{traf_hw}) to our periodic solutions.

If $\hat{a}$ and $n/N$ are specified, from the curves shown in Figure \ref{Ch_m_plot}, the wave modulus $m$ (as well as $a,b,c,d$) and wave speed $\omega$ are identified. Hence, the solution parameters of (\ref{u0jt}) are defined by choosing $\hat{a}$ and $n/N$. Moreover, Figure \ref{Ch_m_plot} reveals that for some fixed $\hat{a}$ and $n/N$, there are two possible values for $\kappa_1/\gamma$, and therefore, two valid fixed point solutions. For the remainder of this paper, solutions with $\kappa_1/\gamma\lesssim 0.128$ ($\kappa_1/\gamma\gtrsim 0.128$) are referred to as the first (second) solution. The first solution represents traffic congestion since it is of a downward form with $u_0\le 1/\sqrt{3}$, which means the headway is less than or equal to the critical headway. Whereas, the second solution is of an upward form with $u_0\ge 1/\sqrt{3}$. Note that when $m\rightarrow 0$, there is only one possible solution ($u_0\rightarrow 1/\sqrt{3}$), where the headway tends to the constant state $h_c$.

Thus, a large family of spatially periodic steady travelling waves have been highlighted that do not slowly evolve, where their amplitude, mean height and period remain fixed. This solution type was discussed by \citet{kom95}, where it was conjectured that they were always unstable. Note however that a different parameter space was considered by \citet{kom95}. Next, the stability of our periodic solutions is examined.
\section{Stability of periodic solutions\label{sps}}
The asymptotic spatially periodic headway solutions are now fully defined, such that
\beq
\Delta x_j(t)=h_c+\epsilon \sqrt{\frac{\omega g_1}{g_2}}\left(u_0(j,t)-\frac{1}{\sqrt{3}}\right),\label{hdwyfull}\eeq
where $u_0(j,t)$ is given by (\ref{u0jt}). The solution parameters $m,a,b,c,d,\omega$ and $\epsilon=\sqrt{(\hat{a}_c/\hat{a})-1}$ are established using the steps outlined in Section \ref{fp} (also refer to Figure \ref{Ch_m_plot}). As previously explained, choosing $\hat{a}$ and $n/N$ then determines the remaining parameters.

\begin{figure}
\begin{center}
\includegraphics[width=4.7in,height=2.5in]{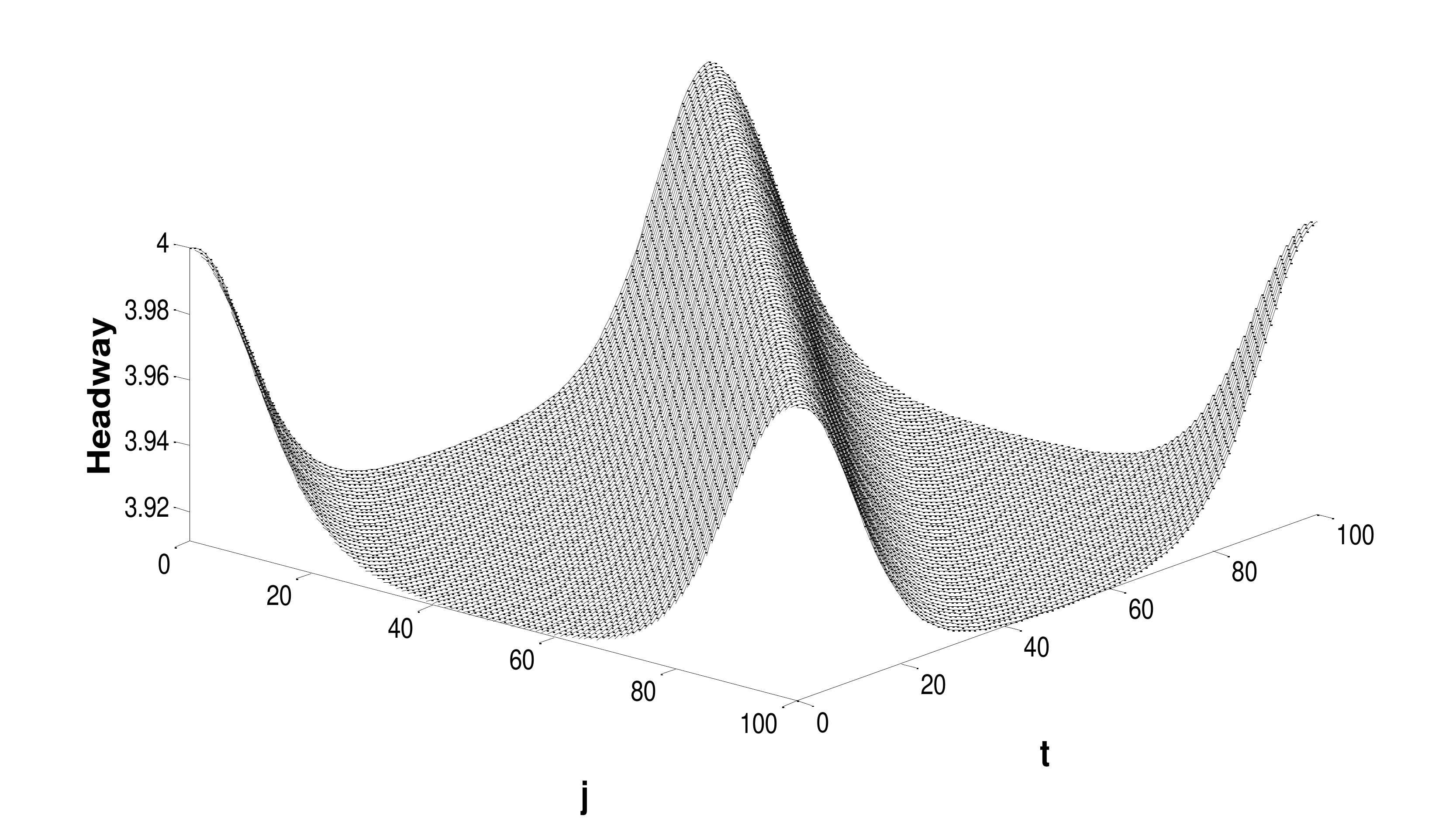}\\
\includegraphics[width=2.8in]{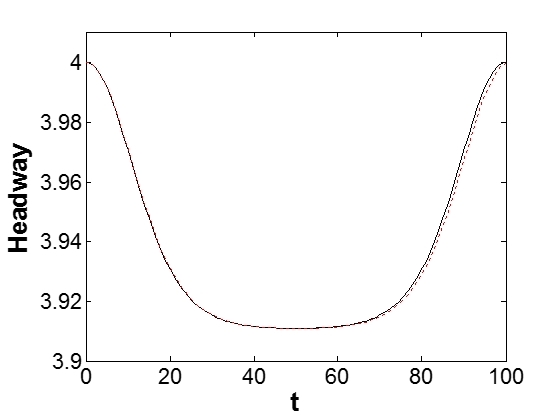}
\includegraphics[width=2.8in]{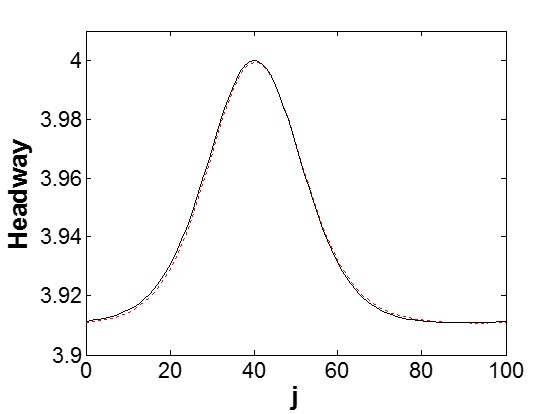}\\
\includegraphics[width=6.6in,height=3.2in]{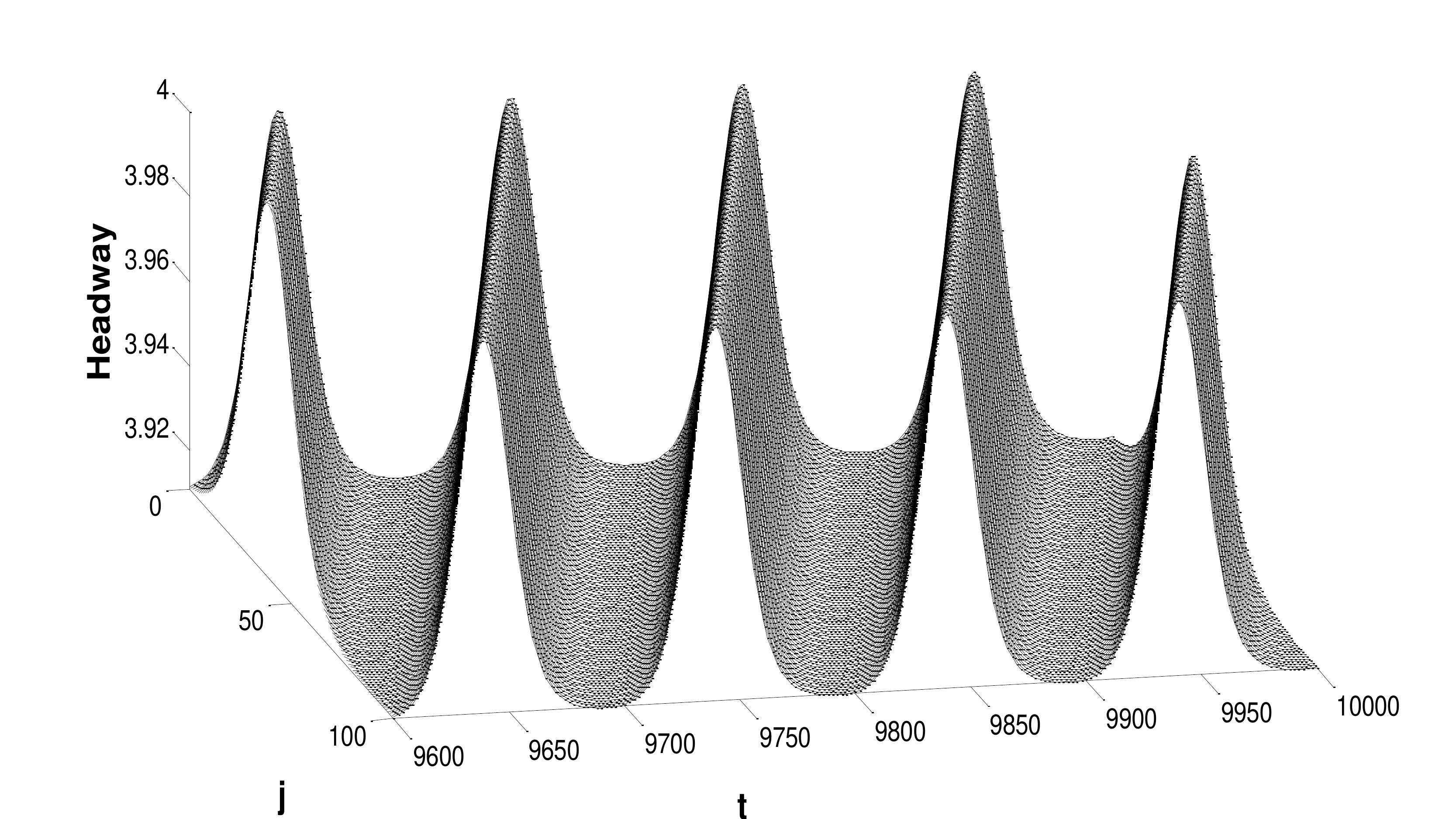}
\caption{Headway solutions to the model (\ref{traf_hw}) that satisfy periodic boundaries, where $m=0.99659$, $\epsilon=0.0709$, $\omega=4.79$, $\hat{a}=1.99$, $n=1$, $N=100$, $\kappa_1/\gamma=0.037582$. Top: The asymptotic solution for $j\in[0,100]$ and $t\in[0,100]$. Middle left: The asymptotic solution (black) compared to the numerical solution (red dotted) when $j=0,100$ and $t\in[0,100]$.  Middle right: The asymptotic solution (black) compared to the numerical solution (red dotted) when $t=10000$ and $j\in[0,100]$. Bottom: The numerical solution for $j\in[0,100]$ and $t\in[9600,10000]$.
\label{eps007n1}}
\end{center}
\end{figure}

\begin{figure}
\begin{center}
\includegraphics[width=4.7in,height=2.5in]{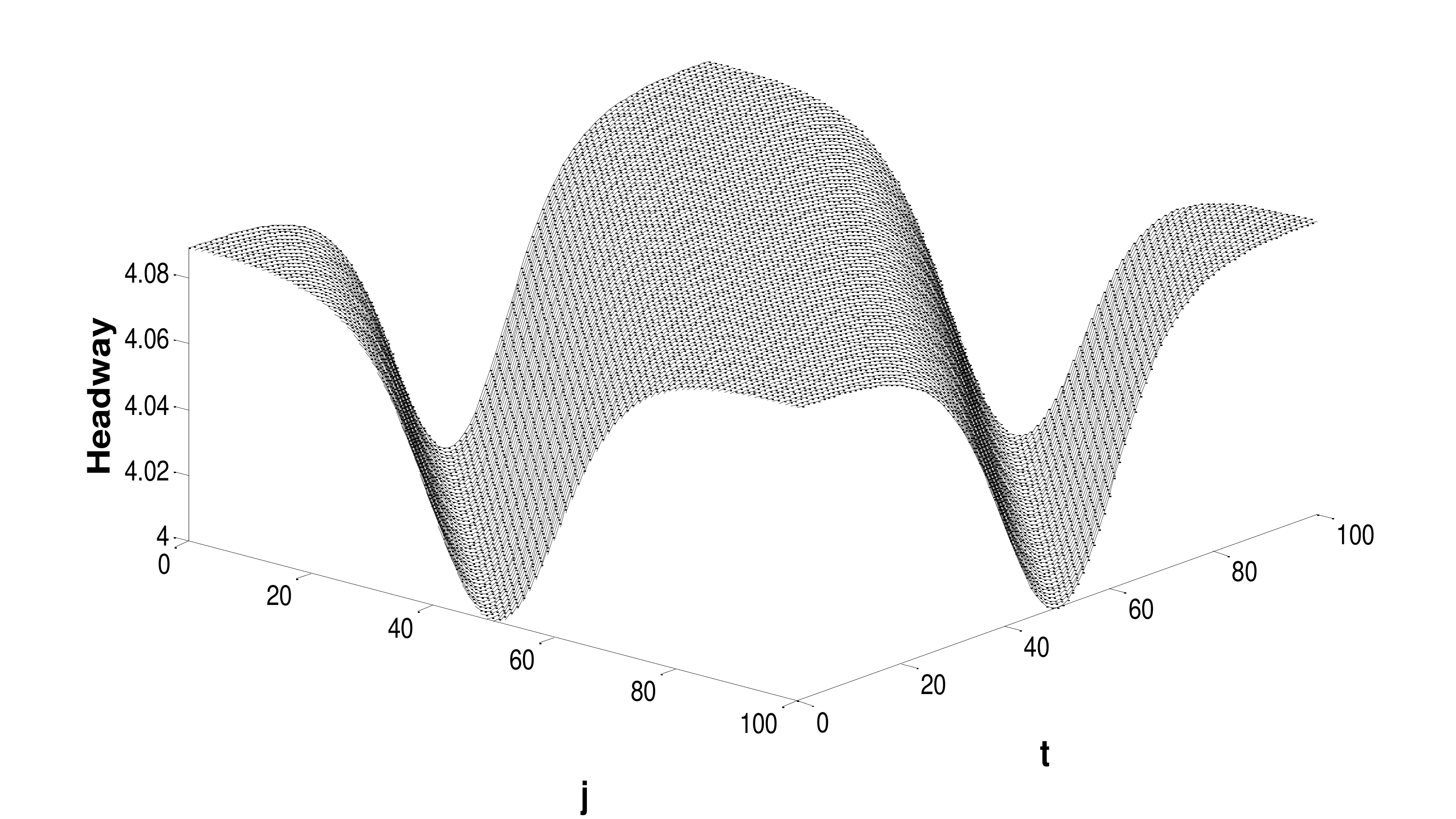}\\
\includegraphics[width=2.8in]{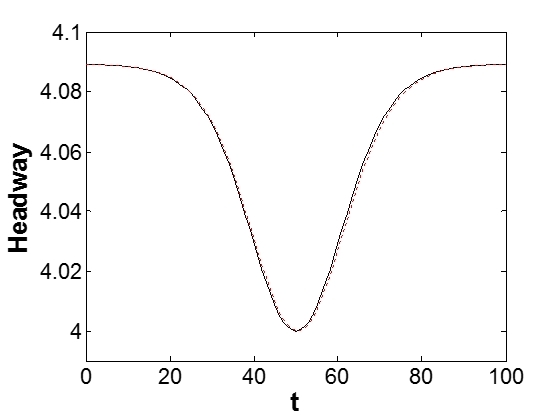}
\includegraphics[width=2.8in]{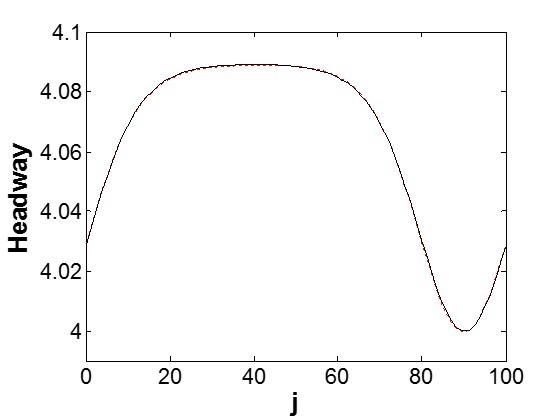}\\
\includegraphics[width=6.5in,height=3.2in]{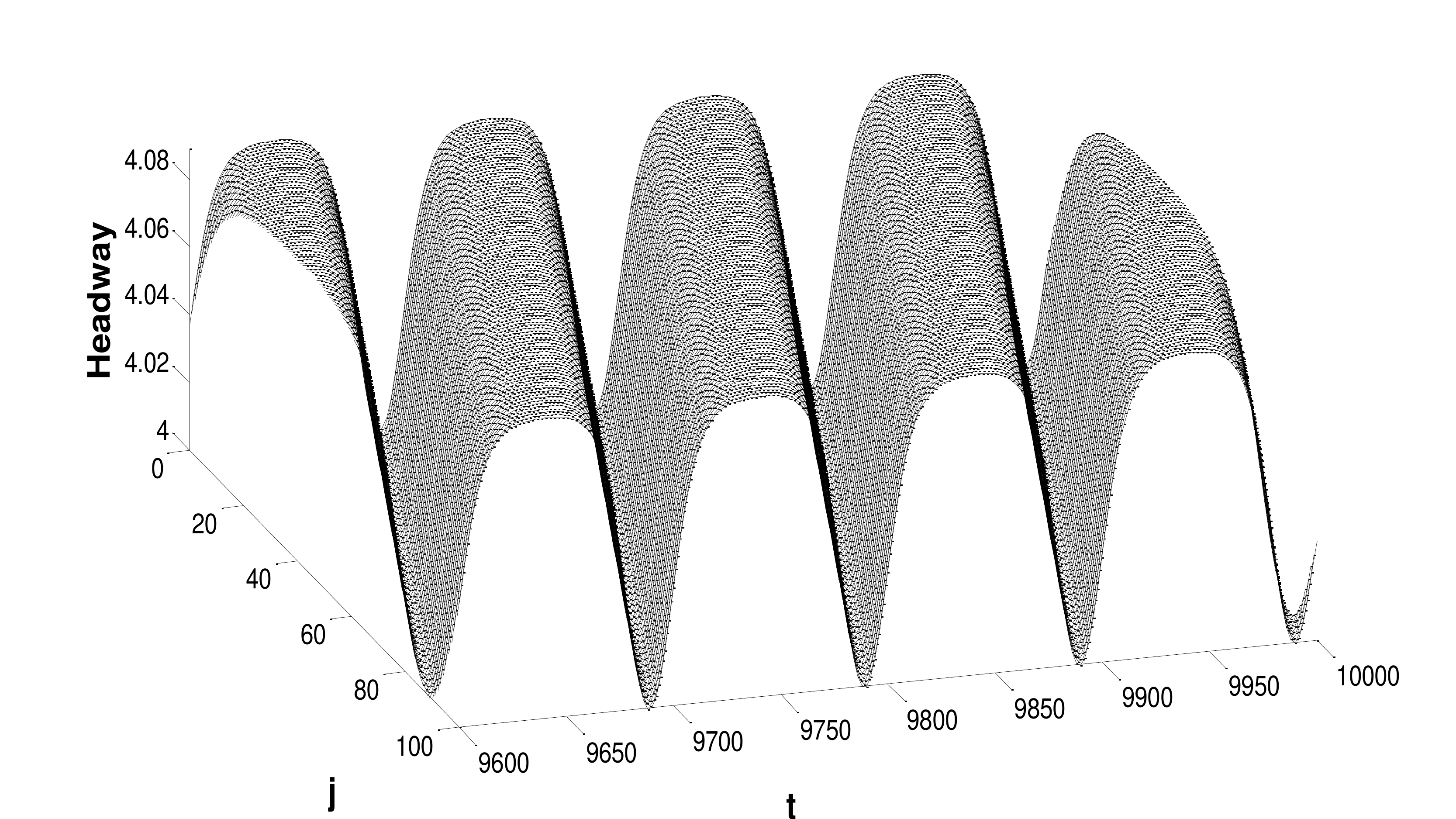}
\caption{Headway solutions to the model (\ref{traf_hw}) that satisfy periodic boundaries, where $m=0.99659$, $\epsilon=0.0709$, $\omega=4.79$, $\hat{a}=1.99$, $n=1$, $N=100$, $\kappa_1/\gamma=0.219018$. Top: The asymptotic solution for $j\in[0,100]$ and $t\in[0,100]$. Middle left: The asymptotic solution (black) compared to the numerical solution (red dotted) when $j=0,100$ and $t\in[0,100]$.  Middle right: The asymptotic solution (black) compared to the numerical solution (red dotted) when $t=10000$ and $j\in[0,100]$. Bottom: The numerical solution for $j\in[0,100]$ and $t\in[9600,10000]$.
\label{eps007n1u}}
\end{center}
\end{figure}

\begin{figure}
\begin{center}
\includegraphics[width=4.7in,height=2.5in]{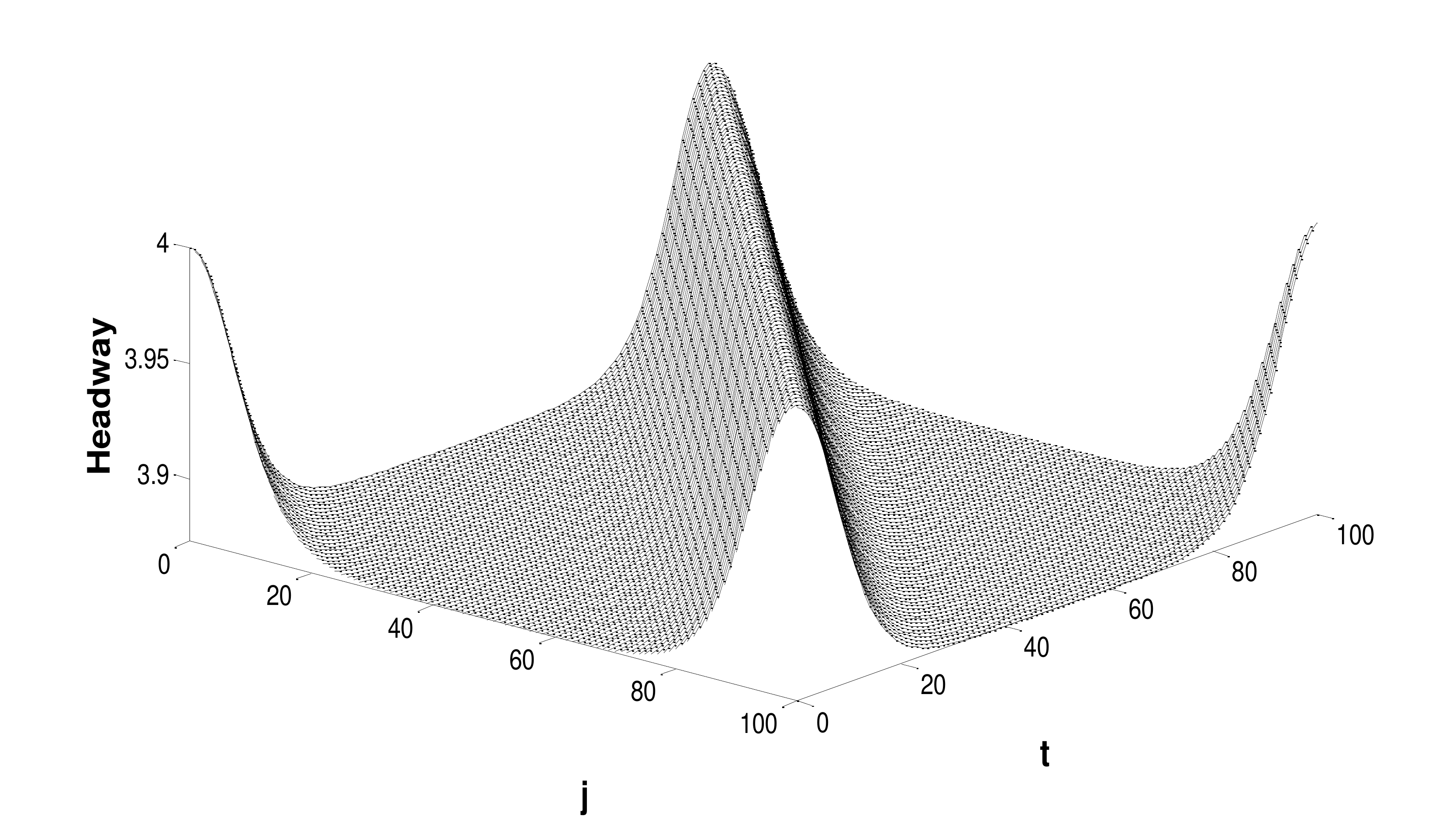}
\includegraphics[width=2.8in]{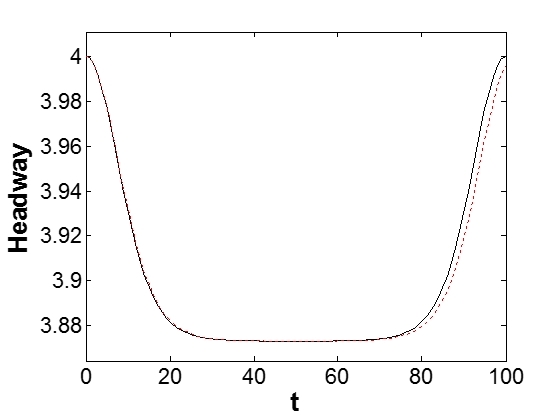}
\includegraphics[width=2.8in]{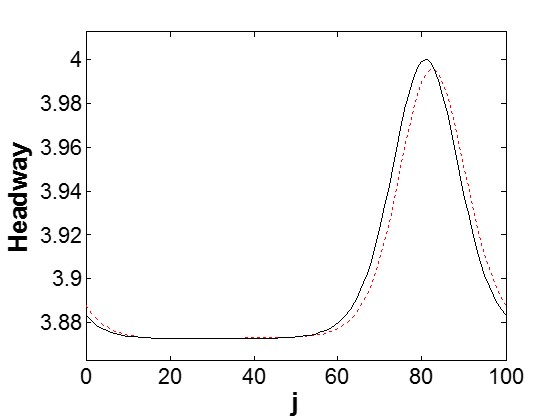}\\
\includegraphics[width=6.5in,height=3.2in]{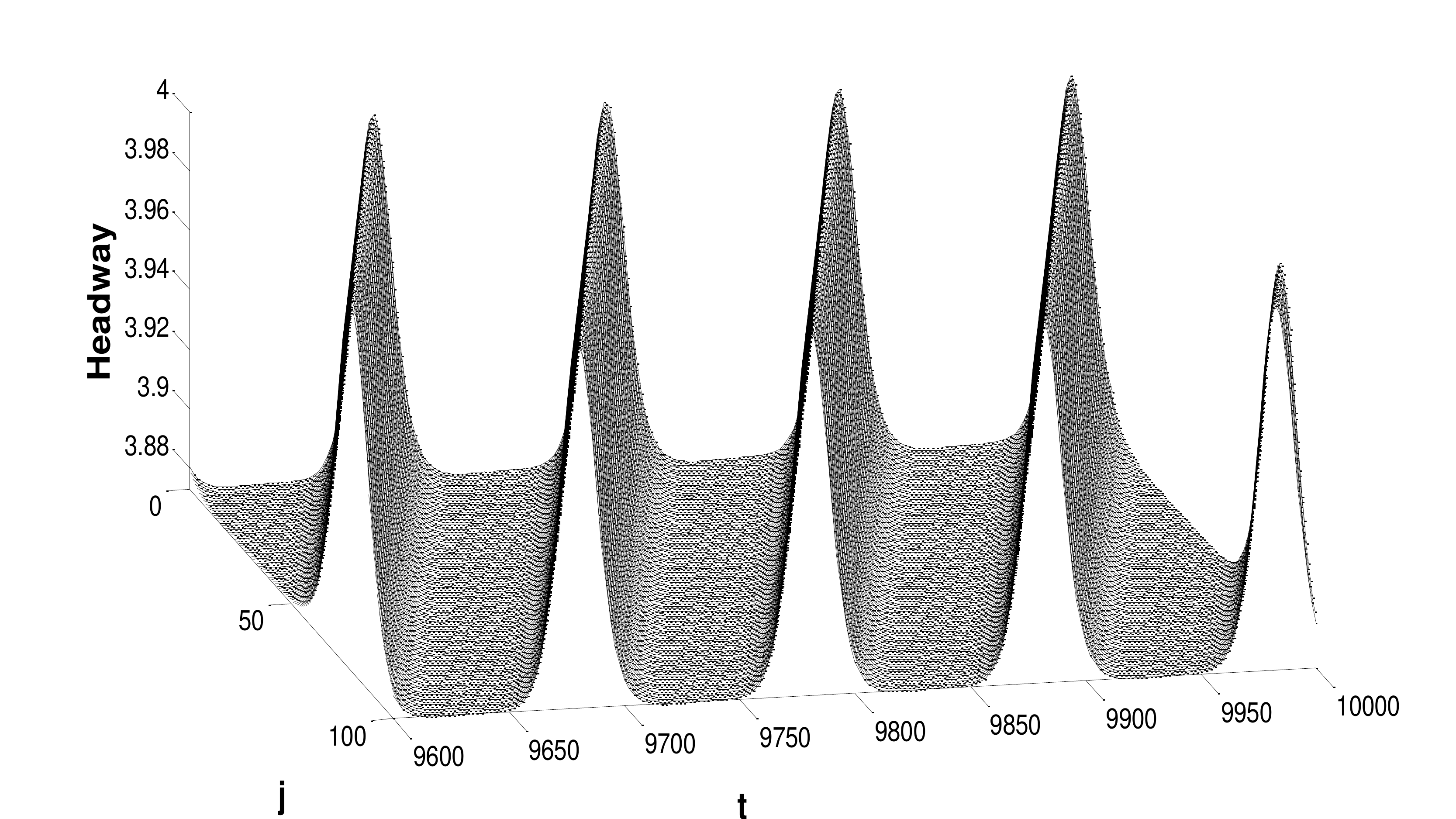}
\caption{Headway solutions to the model (\ref{traf_hw}) that satisfy periodic boundaries, where $m=0.99987$, $\epsilon=0.1005$, $\omega=4.80$, $\hat{a}=1.98$, $n=1$, $N=100$, $\kappa_1/\gamma=0.037578$. Top: The asymptotic solution for $j\in[0,100]$ and $t\in[0,100]$. Middle left: The asymptotic solution (black) compared to the numerical solution (red dotted) when $j=0,100$ and $t\in[0,100]$.  Middle right: The asymptotic solution (black) compared to the numerical solution (red dotted) when $t=10000$ and $j\in[0,100]$. Bottom: The numerical solution for $j\in[0,100]$ and $t\in[9600,10000]$.
\label{eps01n1}}
\end{center}
\end{figure}

\begin{figure}
\begin{center}
\includegraphics[width=3.8in,height=2.5in]{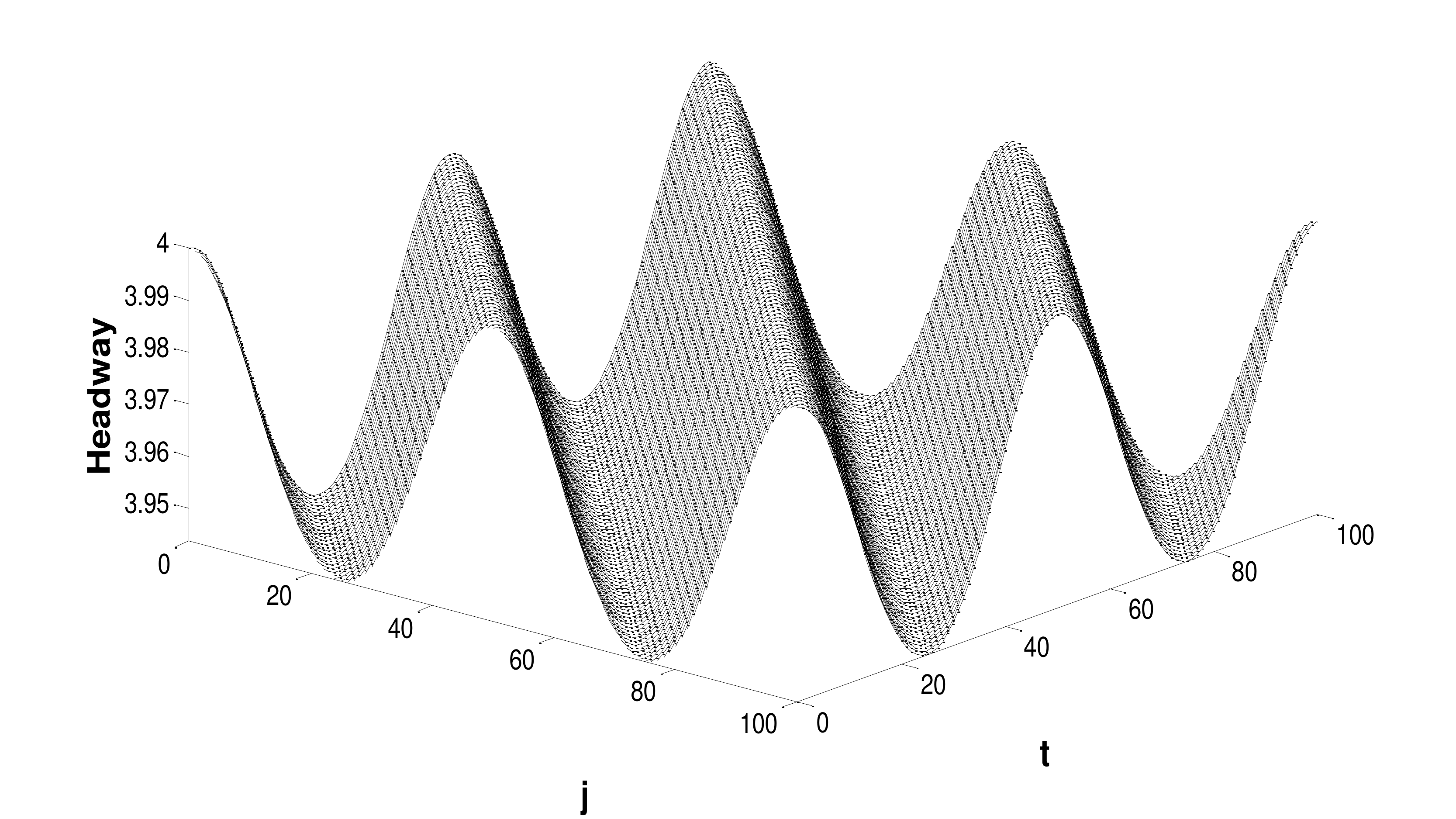}
\includegraphics[width=2.5in]{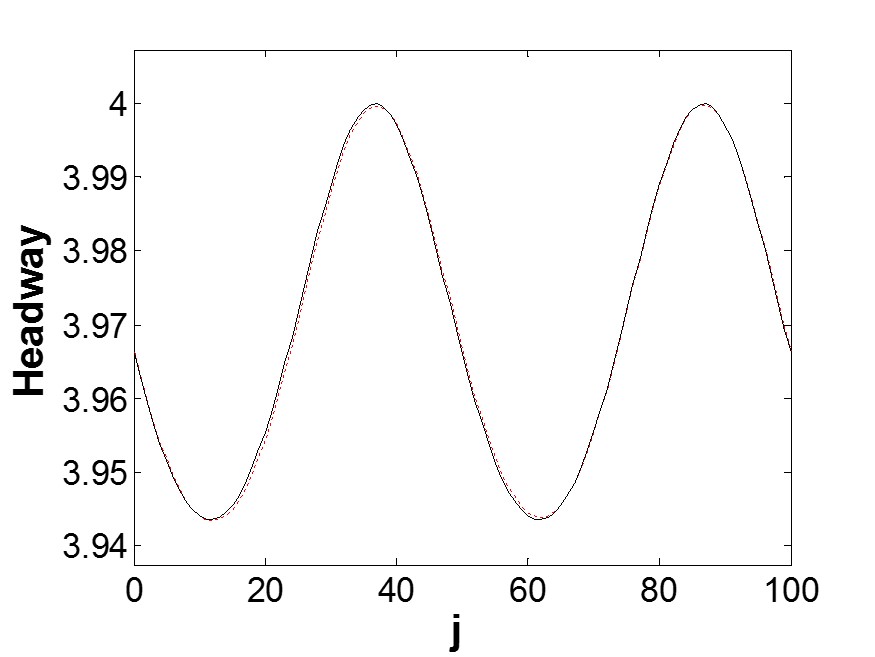}\\
\caption{Headway solutions to the model (\ref{traf_hw}) that satisfy periodic boundaries, where $m=0.792877$, $\epsilon=0.0709$, $\omega=4.38$, $\hat{a}=1.99$, $n=2$, $N=100$, $\kappa_1/\gamma=0.051638$. Left: The asymptotic solution for $j\in[0,100]$ and $t\in[0,100]$. Right: The asymptotic solution (black) compared to the numerical solution (red dotted) when $t=10000$ and $j\in[0,100]$.
\label{eps007n2}}
\end{center}
\end{figure}

The traffic model (\ref{traf_hw}) governing headway is next solved with MATLAB's ode45, where the initial condition is defined by (\ref{hdwyfull}).
As well, periodic boundary conditions are implemented. Then, over different time intervals, the asymptotic and numerical results are compared in Figures \ref{eps007n1}-\ref{eps007n2}. Here, the number of cars on the road is set to $100$.

For Figures \ref{eps007n1}-\ref{eps01n1}, the top panel conveys the asymptotic solution (\ref{u0jt}) over the domain $j\in[0,100]$ and $t\in[0,100]$. The middle panel then compares the asymptotic solution in black to the red-dotted numerical result for $j=0,100$, $t\in[0,100]$ (left) and $j\in[0,100]$, $t=10,000$ (right). Lastly, the bottom plot depicts the numerical solution for $j\in[0,100]$ and $t\in[9600,10000]$.

The first, downward solution for $v_{max}=2$, $h_c=4$, $\hat{a}=1.99$ ($\epsilon=0.0709$), $n=1$ and $N=100$ is shown in Figure \ref{eps007n1}, where $\kappa_1/\gamma=0.037582$. The top panel reveals two distinct zones. These are, a cluster of vehicles with headway $h_{min}<h_c$ represented by the wave trough and a much smaller vehicle cluster with headway $h_c$ corresponding to the wave peak. So, a vehicle travelling at the safety distance $h_c$ decelerates due to a slower preceding vehicle and then endures a prolonged slow period with headway $h_{min}$. They then return to the safety distance momentarily to repeat this process. The middle panel of Figure \ref{eps007n1} reveals excellent agreement between the two solutions. The numerical result at $t$ very large is examined in the middle and lower panels, where the numerical wave appears to propagate without divergence, a change in amplitude or the development of a phase shift, when compared to the asymptotic solution. This suggests that this spatially periodic solution is stable.

Figure \ref{eps007n1u} portrays the alternate second solution for $v_{max}=2$, $h_c=4$, $\hat{a}=1.99$ ($\epsilon=0.0709$), $n=1$ and $N=100$, where $\kappa_1/\gamma=0.219018$. The top panel shows it is of an upward form since now the headway varies between $h_c$ and $h_{max}>h_c$. Therefore, a vehicle travelling at the safety distance $h_c$ will now accelerate as the preceding car is faster. This vehicle then experiences an extended faster state, travelling with the headway $h_{max}$. Next, they decelerate and revert to the safety distance $h_c$. The car remains briefly at $h_c$ and then repeats this motion. The middle panel shows again little discrepancy between the asymptotic and numerical solutions. As well, the middle right and bottom figures investigate the long time behaviour of the numerical result. Once again, no evidence of a phase shift or amplitude variation is exhibited and thus, this solution appears stable.

Next, the driver's sensitivity is reduced to $\hat{a}=1.98$ and as a result, the perturbation parameter increases to $\epsilon=0.1005$. The first solution is displayed in Figure \ref{eps01n1} where $v_{max}=2$, $h_c=4$, $n=1$, $N=100$ and $\kappa_1/\gamma=0.037578$. The vehicle behaviour is consistent with Figure \ref{eps007n1}, except now the wave trough is flatter and $ h_{min}$ has decreased. Hence, the duration of a car travelling with headway $h_{min}$ is longer and the reduction in speed is greater. The middle panel depicts some very small differences between the asymptotic and numerical solutions due to increasing $\epsilon$, although the match is still very good. Analysing the solution at $t$ very large in the middle right and bottom figures, it is apparent that the numerical result is again stable since the phase and amplitude appear constant with increasing $t$. There is also a second, upward stable solution that corresponds to the parameters $\hat{a}=1.98$, $v_{max}=2$, $h_c=4$, $n=1$ and $N=100$, although it is not depicted here. This solution will have the same behaviour as that conveyed by Figure \ref{eps007n1u}, such that $\Delta x_j\in[h_c,h_{max}]$.

So far, the solutions considered all have one peak/trough over the domain $j\in[0,100]$, since $n=1$. Instead, choosing $n>1$, multiple oscillations over the domain will occur. As an example, the downward solution corresponding to $\hat{a}=1.99$, $v_{max}=2$, $h_c=4$, $n=2$, $N=100$, $\kappa_1/\gamma=0.051638$ is shown in Figure \ref{eps007n2}. The behaviour observed on the left is consistent with $n=1$ (Figure \ref{eps007n1}), except waves with two headway troughs/peaks over $j\in[0,100]$ now propagate. The long time dynamics are examined on the right, where the black and red-dotted curves correspond to the asymptotic and numerical solutions respectively when $t=10000$. The solution appears stable since no phase or amplitude changes are exhibited.

Thus, our set of periodic solutions have been shown to be stable using numerical results. In contrast, \citet{kom95} found solutions of a similar form that they supposed were always unstable, and therefore not observed numerically. Note however that our solutions are numerically stable only within a certain neighbourhood of the neutral stability line's critical point (see Section \ref{tfm}).

\section{Conclusion}
The OV model (\ref{traf_hw}) was used to predict traffic behaviour. In particular, this model's linearly unstable region was studied, where (\ref{traf_hw}) transformed into the mKdV equation with higher order correction terms. A multi-scale perturbation method was applied to this equation to locate steady travelling wave solutions. Consequently, the leading order solution that varied with slow and fast variables was defined. A system of differential equations at the next order was also found, which described this solution's slow evolution. The critical points of this system were shown to represent a family of steady travelling waves that had constant amplitude, mean height and period. Imposing periodic boundary constraints then determined the relationship between the solution parameters and the driver's sensitivity, $\hat{a}$, where for some fixed value of $\hat{a}$, two solutions existed of upward and downward form. As a result of establishing this relationship for $\hat{a}$, a numerical investigation was performed. This validated our analysis by demonstrating excellent agreement between the asymptotic and numerical results. Furthermore, we examined the behaviour of these solutions when $t$ was very large. The numerical wave did not diverge, exhibit any phase shift or variation in amplitude, suggesting our set of solutions was stable. Identifying these stable solutions will have important implications for future studies of the traffic model (\ref{traf_hw}), particularly when numerically evaluating (\ref{traf_hw}) within the unstable region and interpreting the results.

\section*{Appendix\label{app}}
The integrals $\alpha_1$ and $\alpha_2$ written in full are
\bsub
\bal
&\alpha_1=\frac{1}{2P}\int_{\theta_1}^{\theta_2}u_0d\theta=c+\frac{(d-c)}{2P}\int_{\theta_1}^{\theta_2}\frac{\operatorname{sn}^2(\beta(\theta-\theta_0))}{e+\operatorname{sn}^2(\beta(\theta-\theta_0))}d\theta,\\
&\alpha_2=\frac{1}{2P}\int_{\theta_1}^{\theta_2}u_0^2d\theta=c^2+\frac{2c(d-c)}{2P}\int_{\theta_1}^{\theta_2}\frac{\operatorname{sn}^2(\beta(\theta-\theta_0))}{e+\operatorname{sn}^2(\beta(\theta-\theta_0))}d\theta
+\frac{(d-c)^2}{2P}\int_{\theta_1}^{\theta_2}\frac{\operatorname{sn}^4(\beta(\theta-\theta_0))}{\left(e+\operatorname{sn}^2(\beta(\theta-\theta_0))\right)^2}d\theta,
\end{align*}\esub
given (\ref{u0soln}). Using \citet{byr54} to solve these integrals, we obtain
\bsub
\bal
&\frac{1}{2P}\int_{\theta_1}^{\theta_2}\frac{\operatorname{sn}^2(\beta(\theta-\theta_0))}{e+\operatorname{sn}^2(\beta(\theta-\theta_0))}d\theta =1-\frac{\Pi(-1/e,m)}{K(m)},\\
&\frac{1}{2P}\int_{\theta_1}^{\theta_2}\frac{\operatorname{sn}^4(\beta(\theta-\theta_0))}{\left(e+\operatorname{sn}^2(\beta(\theta-\theta_0))\right)^2}d\theta \\ &=\frac{1}{K(m)}\left(K(m)-2\Pi(-1/e,m)\right)\\
&+\frac{1}{2K(m)(-1/e-1)(m^2+1/e)}\left(-\frac{E(m)}{e}+(m^2+1/e)K(m)+\Pi(-1/e,m)(-2m^2/e-2/e-1/e^2-3m^2)\right),
\end{align*}\esub
where $K$, $E$ and $\Pi$ are the complete elliptic integrals of the first, second and third kind respectively. Therefore, $\alpha_1$ and $\alpha_2$ are functions of the solution parameters $m$, $c$, $d$ and $e$ (see (\ref{u0soln})).
\bibliographystyle{plainnat}
\bibliography{bib_lhattam}

\begin{thebibliography}{12}
\providecommand{\natexlab}[1]{#1}
\providecommand{\url}[1]{\texttt{#1}}
\expandafter\ifx\csname urlstyle\endcsname\relax
  \providecommand{\doi}[1]{doi: #1}\else
  \providecommand{\doi}{doi: \begingroup \urlstyle{rm}\Url}\fi

\bibitem[Bando et~al.(1995)Bando, Hasebe, Nakayama, Shibata, and
  Sugiyama]{ban95}
M.~Bando, K.~Hasebe, A.~Nakayama, A.~Shibata, and Y.~Sugiyama.
\newblock Dynamical model of traffic congestion and numerical simulation.
\newblock \emph{Phys. Rev. E}, 51:\penalty0 1035--1042, 1995.

\bibitem[Byrd and Friedman(1954)]{byr54}
P.~Byrd and M.~Friedman.
\newblock \emph{Handbook of Elliptic Integrals for Engineers and Physicists}.
\newblock Springer-Verlag, Berlin, 1954.

\bibitem[Ge et~al.(2005)Ge, Cheng, and Dai]{ge05}
H.X. Ge, R.J. Cheng, and S.Q. Dai.
\newblock {K}d{V} and kink-antikink solitons in car-following models.
\newblock \emph{Physica A}, 357\penalty0 (3–4):\penalty0 466--476, 2005.

\bibitem[Hattam(2016)]{hat16}
L.L. Hattam.
\newblock {K}d{V} cnoidal waves in a traffic flow model with periodic
  boundaries.
\newblock \emph{Physica D (under review)}, 2016.

\bibitem[Hattam and Clarke(2015)]{hat15}
L.L. Hattam and S.R. Clarke.
\newblock Modulation theory for the steady forced {K}d{V}-{B}urgers equation
  and the construction of periodic solutions.
\newblock \emph{Wave Motion}, 56\penalty0 (0):\penalty0 67--84, 2015.

\bibitem[Kamchatnov et~al.(2012)Kamchatnov, Kuo, Lin, Horng, Gou, Clift, El,
  and Grimshaw]{kam12}
A.M. Kamchatnov, Y.-H. Kuo, T.-C. Lin, T.-L. Horng, S.-C. Gou, R.~Clift, G.A.
  El, and R.H.J. Grimshaw.
\newblock Undular bore theory for the {G}ardner equation.
\newblock \emph{Phys. Rev. E}, 86:\penalty0 036605--036627, 2012.

\bibitem[Komatsu and Sasa(1995)]{kom95}
T.S. Komatsu and S.~Sasa.
\newblock Kink soliton characterizing traffic congestion.
\newblock \emph{Phys. Rev. E}, 52:\penalty0 5574--5582, 1995.

\bibitem[Li et~al.(2015)Li, Zhang, Xu, and Qian]{li15}
Z.~Li, R.~Zhang, S.~Xu, and Y.~Qian.
\newblock Study on the effects of driver's lane-changing aggressiveness on
  traffic stability from an extended two-lane lattice model.
\newblock \emph{Communications in Nonlinear Science and Numerical Simulation},
  24\penalty0 (1-3):\penalty0 52--63, 2015.

\bibitem[Muramatsu and Nagatani(1999)]{mur99}
M.~Muramatsu and T.~Nagatani.
\newblock Soliton and kink jams in traffic flow with open boundaries.
\newblock \emph{Phys. Rev. E}, 60:\penalty0 180--187, 1999.

\bibitem[Nagatani(2002)]{nag02}
T.~Nagatani.
\newblock The physics of traffic jams.
\newblock \emph{Rep. Prog. Phys.}, 65\penalty0 (9):\penalty0 1331, 2002.

\bibitem[Zheng et~al.(2012)Zheng, Tian, Sun, and Liu]{zhe12}
L.-J. Zheng, C.~Tian, D.-H. Sun, and W.-N. Liu.
\newblock A new car-following model with consideration of anticipation driving
  behavior.
\newblock \emph{Nonlinear Dynamics}, 70\penalty0 (2):\penalty0 1205--1211,
  2012.

\bibitem[Zhu and Dai(2008)]{zhu08}
H.B. Zhu and S.Q. Dai.
\newblock Numerical simulation of soliton and kink density waves in traffic
  flow with periodic boundaries.
\newblock \emph{Physica A}, 387\penalty0 (16-17):\penalty0 4367--4375, 2008.

\end{thebibliography}

\end{document}